\documentclass[12pt]{article}
\usepackage{amsfonts}
\usepackage{amssymb,amsfonts,amsmath, psfrag,eepic,colordvi,epsfig}
by \topskip 
\numberwithin{equation}{section}
\newtheorem{theorem}{Theorem}[section]

\newtheorem{lemma}[theorem]{Lemma}

\begin{document}
\parskip 6pt

\pagenumbering{arabic}
\def\sof{\hfill\rule{2mm}{2mm}}
\def\ls{\leq}
\def\gs{\geq}
\def\SS{\mathcal S}
\def\qq{{\bold q}}
\def\MM{\mathcal M}
\def\TT{\mathcal T}
\def\EE{\mathcal E}
\def\lsp{\mbox{lsp}}
\def\rsp{\mbox{rsp}}
\def\pf{\noindent {\it Proof.} }
\def\mp{\mbox{pyramid}}
\def\mb{\mbox{block}}
\def\mc{\mbox{cross}}
\def\qed{\hfill \rule{4pt}{7pt}}
\def\block{\hfill \rule{5pt}{5pt}}
\begin{center}
{\large {\bf Faulhaber's Theorem on Power Sums}}

\vskip 6mm

{\small William Y.C. Chen$^1$, Amy M. Fu$^2$, and Iris F. Zhang$^3$\\[%
2mm] Center for Combinatorics, LPMC-TJKLC\\
Nankai University, Tianjin 300071,
P.R. China \\[3mm]
$^1$chen@nankai.edu.cn, $^2$fu@nankai.edu.cn, $^3$zhangfan03@mail.nankai.edu.cn \\[0pt%
] }
\end{center}

\vskip10mm
\noindent {\bf Abstract.}  We observe that the classical Faulhaber's
theorem on sums of odd powers also holds for an arbitrary arithmetic
progression, namely, the odd power sums of any arithmetic
progression $a+b, a+2b, \ldots, a+nb$ is a polynomial in
$na+n(n+1)b/2$. While this assertion can be deduced from the
original Fauhalber's theorem, we give an alternative formula
  in terms of the Bernoulli
polynomials. Moreover, by utilizing the central factorial numbers as
in the approach of Knuth, we derive formulas for $r$-fold sums of
powers without resorting to the notion of $r$-reflexive functions.
We also provide  formulas for the $r$-fold alternating sums of
powers in terms of Euler polynomials.




\vskip10mm

\section{Introduction}

The classical theorem of Faulhaber  states that the sums of odd
powers
\[ 1^{2m-1}+2^{2m-1}+\ldots+n^{2m-1}\]
 can be
expressed as a polynomial of the triangular number $T_n=n(n+1)/2$;
See Beardon \cite{Beardon}, Knuth \cite{ knuth}. Moreover, Faulhaber
observed that the $r$-fold summation of $n^{m}$ is a polynomial in
$n(n+r)$ when $m$ is positive
 and $m-r$ is even \cite{knuth}.
The classical Faulhaber theorem for odd power sums
 was proved by Jacobi \cite{jacobi}; See also Edwards
\cite{edwards}. Let us recall the notation on the $r$-fold power
sums: $\textstyle \sum^0 n^m = n^m$, and
\begin{equation} \label{rfold-d}
{\sum}^{r} n^m = {\sum}^{r-1} 1^m + {\sum}^{r-1} 2^m + \cdots+
{\sum}^{r-1} n^m .
\end{equation}
For example, $\sum^1 n^m = 1^m +2^m +\cdots + n^m$, and
\begin{equation*}
{\textstyle  \sum^2 n^m =  \sum^1 1^m + \sum^1 2^m +\cdots + \sum^1 n^m}=
                 \sum_{i=1}^n (n+1-i) i^m .
\end{equation*}
For even powers, it has been shown that the sum $1^{2m}+2^{2m} +
\cdots + n^{2m}$ is a polynomial in the triangular number $T_n$
multiplied by a linear factor in $n$. Gessel and Viennot
\cite{gessel} had a remarkable discovery that the alternating sum
 ${\sum^n_{i=1}}(-1)^{n-i}i^{2m}$ is also a polynomial in
the triangular number $T_n$.

Faulhaber's theorem has drawn much attention from various points of
view. Grosset and  Veselov \cite{grosset} investigated a
generalization of the Faulhaber polynomials related to elliptic
curves.  Warnaar \cite{warnaar}, Schlosser \cite{schlosser} and Zhao
and Feng \cite{ZF} studied the $q$-analogues of the formulas for the
first few power sums. Garrett \cite{garrett} found a combinatorial
proof of the formula for sums of $q$-cubes. Guo and Zeng \cite{guo}
obtained the $q$-analogue formula in the general case. Furthermore,
Guo, Rubey and Zeng \cite{GRZ} have shown that the $q$-Faulhaber and
$q$-Sali\'e coefficients are nonnegative and symmetric in a
combinatorial setting of nonintersecting lattice paths.

In this paper, we first formulate Faulhaber's theorem in a more
general framework, that is, in terms of power sums of an arithmetic
progression. Given an arithmetic progression:
\[ a+b, \;\; a+2b, \;\; \ldots, \;\; a+nb,\]
Faulhaber's theorem implies that  odd power sums of the above series
are polynomials in $na+n(n+1)b/2$. In particular, an odd power sum
of the first $n$ odd numbers \[ 1^{2m-1}+ 3^{2m-1} + \cdots+
(2n-1)^{2m-1}\]
 is a
polynomial in $n^2$,
 and the sum
\[ 1^{2m-1} + 4^{2m-1} + 7^{2m-1} +\cdots + (3n-2)^{2m-1}\]
is a polynomial in the pentagonal number $n(3n-1)/2$.

Because of the relation $(a+bi)^m = b^m (a/b + i)^m$, there is no
loss of generality to consider the series \begin{equation}\label{xn}
 x+1, x+2,
\ldots, x+n.\end{equation}
 Let
 \begin{equation}\label{ln}
 \lambda=n(n+2x+1)
 \end{equation}
 be the sum of the
sequence $x+1, x+2, \ldots, x+n$. The the power sums
\[ S_{2m-1}
=(x+1)^{2m-1} +(x+2)^{2m-1} + \cdots + (x+n)^{2m-1}\] is a
polynomial in $\lambda$. For exmaple,
\begin{eqnarray*}
S_3&=& {\lambda^2\over 4}+{(x^2+x)\over 2}\lambda;\\
S_5&=&{\lambda^3\over6} +{1\over12}(6x^2 +6x - 1)\lambda^2
       +{1\over6}(3x^4+6x^3 +2x^2 -x )\lambda;\\
S_7&=& {\lambda^4\over 8}+{1\over 6}(3x^2+3x   -1)\lambda^3+
{1\over 12}(9x^4 + 18x^3+ 3x^2   -6x+1 )\lambda^2\\
& &+\,{1\over 6}(3x^6+9x^5 +6x^4 -3x^3 -2x^2 +x )\lambda.
\end{eqnarray*}

 It
should be noticed that the above more general setting of Faulhaber's
theorem can be deduced on the the original version of Faulhaber's
theorem. When $x$ is a positive integer, we have the relation
\[
\sum_{i=1}^{n} (x+i)^m= \sum^{n+x}_{i=1}i^m-\sum^{x}_{i=1}i^m.
\]
By Faulhaber's theorem,
 $\sum^{n+x}_{i=1}i^m$ and
$\sum^{x}_{i=1}i^m$ are polynomials in $(n+x)(n+x+1)$ and $x(x+1)$,
respectively. Using the following simple but important relation
 \begin{equation}
  (n+x)(n+x+1)=n(n+2x+1) +
x(x+1),\end{equation}
 we see that
\begin{equation}\label{lam}
[(n+x)(n+x+1)]^i-[x(x+1)]^i=\sum^i_{k=1}{i\choose
k}[n(n+2x+1)]^k[x(x+1)]^{i-k},
\end{equation}
which is a polynomial in $n(n+2x+1)$. Clearly, one sees that the
above assertion holds for all real numbers $x$.

Although in principle Faulhaber's theorem is valid for any
arithmetic progression, from a computational point of view it still
seems worthwhile to find a formula for the coefficients in terms
Bernoulli polynomials. The main result of this paper is an approach
to the study of the $r$-fold sums of powers without resorting to the
properties of $r$-reflective functions as in the approach of Knuth
\cite{knuth}. In the last section, we obtain  formulas for the
$r$-fold alternating sums of powers in terms of the Euler
polynomials.

\section{An Alternative Formula}

In this section, we give an explicit formula for the coefficients
regarding Faulhaber's theorem for the series $x+1, x+2, \ldots,
x+n$, which reduces to an alternative formula to the Gessel-Viennot
formula when setting $x=0$. We first recall some basic facts about
Bernoulli polynomials $B_n(x)$ which are defined by the following
generating function:
\begin{equation}\label{bgf}
\sum^\infty_{n=0}\frac{B_n(x)t^n}{n!}=\frac{te^{xt}}{e^t-1}.
\end{equation}
The power sums of the first $n$ positive integers can be expressed
in terms of $B_i(x)$:
$$\sum_{i=1}^n i^m=\frac{1}{m+1}(B_{m+1}(n+1)-B_{m+1}(1)).$$
Moreover, we have
\begin{equation}\label{bnx}
\sum_{i=1}^n (x+i)^{2m-1}={1\over 2m}\left( B_{2m}(x+n+1)
-B_{2m}(x+1)\right).
\end{equation}

The Bernoulli numbers $B_n$ are  given by $B_n=B_n(0)$. Note that
the Bernoulli polynomials satisfy the following relations
\begin{eqnarray}\label{eqn52}
B_n(x+1)-B_n(x)&= &nx^{n-1},\\[6pt]
\label{1xx} B_n(1-x)& =& (-1)^nB_n(x), \\[6pt]
\label{de} \frac{d}{dx}B_n(x)& =& nB_{n-1}(x),\\[6pt]
 \label{eqn51}
B_n(x+y)& =& \sum^n_{i=0}{n\choose i}B_i(x)y^{n-i}.
\end{eqnarray}
The evaluation of Bernoulli polynomials at $1/2$ is of special
interest. For $n\geq 0$, we have
\begin{equation}\label{b2n}
B_{2n+1}(1/2)=0, \qquad B_{2n}(1/2)=(2^{1-2n}-1)B_{2n} .
\end{equation}
 From (\ref{eqn51}) and (\ref{b2n}), we can deduce the following
 form of Faulhaber's theorem.

\begin{theorem}
 Let $\lambda=n(n+2x+1)$. Then we have
\begin{eqnarray}\label{sum(x+n)^2m-1}
\sum_{i=1}^n (x+i)^{2m-1}=\sum^m_{k=1}F^{(m)}_k(x) \lambda^k,
\end{eqnarray}
where
\begin{equation}\label{am}
F^{(m)}_k(x)=\frac{1}{2m}\sum^m_{i=k} {2m\choose 2i}{i\choose
k}\left(x+{1\over 2}\right)^{2i-2k} B_{2m-2i}\left(1\over2\right).
\end{equation}
\end{theorem}

\pf From the binomial expansion (\ref{eqn51}), we get
\begin{equation}
B_{2m}(x+n+1)=\sum^{2m}_{i=0}{2m\choose
i}B_{2m-i}\left(1\over2\right)\cdot\left(x+n+{1\over2}\right)^i.
\end{equation}
It follows from (\ref{bnx}) and  (\ref{b2n}) that
\begin{eqnarray*}
\sum(x+n)^{2m-1}={1\over 2m}\sum^{m}_{i=0}{2m\choose
2i}B_{2m-2i}\left(1\over2\right)\cdot\left(x+n+{1\over2}\right)^{2i}-
{1\over 2m} B_{2m}(x+1).
\end{eqnarray*}
Since
\begin{eqnarray}
\left(x+n+{1\over2}\right)^{2i}=\left(\lambda+\left(x+{1\over
2}\right)^2\right)^i=\sum^i_{k=0}{i\choose k}\left(x+{1\over
2}\right)^{2i-2k}\lambda^k,\label{lambda}
\end{eqnarray}
we immediately get (\ref{am}) for $k\geqslant 1$. For $k=0$, we have
\begin{eqnarray*}
F^{(m)}_0(x) ={1\over 2m}\sum^m_{i=0}{2m\choose
2i}B_{2m-2i}\left({1\over
2}\right)\left(x+{1\over2}\right)^{2i}={1\over 2m}B_{2m}(x+1).
\end{eqnarray*}
This completes the proof. \qed

The above formula can be viewed as an alternative form of the
formula of Gessel and Viennot \cite{gessel} for the coefficients
$A_k^{(m)}=F_k^{(m)}(0)$:
\[ A_k^{(m)}= (-1)^{m-k} \sum_j {2m \choose m-k-j}
{m-k+j\choose j} {m-k-j \over m-k+j} B_{m+k+j}, \quad 0\leq k <m.\]

Note that $B_0=1$ and $B_1=-1/2$ are used in the above formula
whereas the formula (\ref{am}) does not involve $B_1$. The
equivalence between the formulas for $F^{(m)}_k(0)$ and $A_k^{(m)}$
can be established via the following generating function for the
coefficients $F_k^{(m+1)}(x)$. The proof is analogous to that given
by Gessel and Viennot \cite{gessel} for the case $x=0$.

\begin{theorem} We have
\begin{equation}\label{fij}
\sum^\infty_{m=0}\sum^\infty_{k=1}F_k^{(m+1)}(x)
t^k\frac{y^{2m+1}}{(2m+1)!}=\frac{\cosh y\sqrt{(x+{1\over
2})^2+t}-\cosh y(x+{1\over 2})}{2\sinh({y\over2})}.
\end{equation}
\end{theorem}


Similarly, the theorem of Gessel and Viennot on  the alternating
sums of even powers  of the first $n$ natural numbers can be
extended to an arithmetic progression $x+1, x+2, \ldots, x+n$. It
turns out that the Euler polynomials play the same role as the
Bernoulli polynomials for sums of odd powers.

The Euler polynomials $E_n(x)$ are defined by
\begin{equation} \label{e-g}
\sum^\infty_{n=0}E_n(x)\frac{t^n}{n!}=\frac{2e^{xt}}{e^t+1}.
\end{equation}
The following expansion formula holds:
\begin{equation}\label{e-bi}
E_n(x+y)=\sum_{k=0}^n{n\choose k}E_k(x)y^{n-k}.
\end{equation}
 For positive even number $n$, we have $E_{n}(1)=0$ . The Euler
numbers $E_n$ and the Euler polynomials are related by
 \begin{equation}\label{e-n-p}
E_{2n+1}=0, \quad E_n=2^nE_n(1/2), \quad n\geq 0.
\end{equation}

\begin{theorem}
Let $\lambda=n(n+2x+1)$. Then we have
\begin{eqnarray*}
\sum^n_{i=1}(-1)^{n-i}(x+i)^{2m}=\sum^m_{k=0}G^{(m)}_k(x)\lambda^k,
\end{eqnarray*}
where
\begin{eqnarray*}
G^{(m)}_k(x)&=&{1\over 2}\sum^m_{i=k}{2m\choose 2i}{i\choose
k}E_{2m-2i}\left({1\over 2}\right)\left(x+{1\over
2}\right)^{i-k},\qquad  1\leq k\leq m;\\
G^{(m)}_0(x)&=&{1\over 2}(1-(-1)^n)E_{2m}(x+1).
\end{eqnarray*}
\end{theorem}

The generating function for $G_k^{(m)}(x)$ is given below, which is
a straightforward extension of the formula of Gessel and Viennot
\cite{gessel}.

\begin{theorem}\label{sgf} We have
\begin{equation} \label{g-ij}
 \sum^\infty_{m=0}\sum^\infty_{k=0}G_k^{(m)}(x)
t^k\frac{y^{2m}}{(2m)!}=\frac{\cosh y\sqrt{(x+{1\over
2})^2+t}}{2\cosh({y\over2})}.
\end{equation}
\end{theorem}

\section{$r$-Fold Sums of Powers}

In this section,  we derive a formula for the $r$-fold sums of
powers of the series $x+1, x+2, \ldots, x+n$ in terms of the central
factorial numbers as used in the approach of Knuth \cite{knuth}. A
key step in our approach is the $r$-fold summation formula for the
lower factorials.  It can be seen from our formula that if $r$ and
$m$ have the same parity then the $r$-fold power sum
${\sum}^r(n+x)^m$ is a polynomial in $n(n+2x+r)$ plus a term that
vanishes when $x=0$. The cases when $r$ and $m$ have different
parities can be dealt with some care, and the details are omitted.

Recall the notation for the lower factorials
$(x)_k=x(x-1)\cdots(x-k+1)$. The definition of the central
factorials $x^{[k]}$ \cite[p. 213]{riordan}, is given by
\[ x^{[k]}=x(x+k/2-1)_{k-1} . \] The central factorial numbers
$T(m, k)$ are determined  by the following relation:
\begin{equation} \label{c-d}
x^m=\sum^m_{k=1}T(m,k)x^{[k]},\quad m\geq 1.
\end{equation}
 Note that $T(m,k)=0$ when $m-k$ is odd. In
particular, we need the following relation
\begin{equation}\label{x^(2m-1)}
x^{2m-1}=\sum^m_{k=1}T(2m,2k)(x+k-1)_{2k-1}.
\end{equation}

We first give a formula for the $r$-fold sums of lower factorials.
From the recursive definition (\ref{rfold-d}) of $r$-fold
summations, we have
\begin{eqnarray} \label{rfold-lf}
{\sum}^r(x+n)_l=\sum_{1\leqslant i_1\leqslant
i_2\leqslant\dots\leqslant i_r\leqslant n}(x+i_1)_l.
\end{eqnarray}

The above multiple summation can be simplified to a single sum.

\begin{theorem}\label{sumrfacto} We have
\begin{eqnarray*}
{\sum}^r(x+n)_{l}&=&\frac{(x+n+r)_{l+r}}{(l+r)_r}-\sum^r_{i=1}{n+r-i-1\choose
r-i}\frac{(x+i)_{l+i}}{(l+i)_i}.
\end{eqnarray*}
\end{theorem}
\pf Setting
 $$F(n)=\frac{1}{l+1}(x+n+1)_{l+1}$$ gives
$F(i)-F(i-1)=(x+i)_l.$ Hence we get
\begin{eqnarray}
\sum_{1\leqslant i_1\leqslant i_2}(x+i_1)_l&=&\sum^{i_2}_{i_1=1}\Big(F(i_1)-F(i_1-1)\Big)\nonumber\\
&=&\frac{1}{l+1}(x+i_2+1)_{l+1}-\frac{1}{l+1}(x+1)_{l+1}.\label{l-f}
\end{eqnarray}
Iterating  (\ref{l-f})  $r-1$ times and using the following
identity
$$\sum^n_{i=1}{l+i-1 \choose l}={l+n \choose l+1},$$  we obtain the desired formula.
\qed

From Theorem \ref{sumrfacto} and the relation (\ref{x^(2m-1)}), we
derive two formulas for the $r$-fold sums of the $m$-th powers when
$r$ and $m$ have the same parity.

\begin{theorem}\label{thm3_1} For $m\geq 1$, we have
\begin{eqnarray*}
{\sum}^{2r+1}(x+n)^{2m-1}&=&\sum^m_{k=1}T(2m,2k)\Bigg\{\frac{(x+n+k+2r)_{2k+2r}}{(2k+2r)_{2r+1}}\\
&&\quad\quad\quad\quad-\sum^{2r+1}_{i=1}{n+2r-i\choose
2r-i+1}\frac{(x+k+i-1)_{2k+i-1}}{(2k+i-1)_{i}}\Bigg\}.
\end{eqnarray*}
\end{theorem}

We remark that  the second summation in the above formula vanishes
when $x=0$, and the lower factorial $(n+k+2r)_{2k+2r}$ can be
rewritten as
\[
 \prod^{k+r}_{i=1}(n+k+2r+1-i)(n-k+i)=
\prod^{k+r}_{i=1}[n(n+2r+1)-(k+2r-i+1)(k-i)]. \] Hence we obtain
Faulhaber's theorem for the $(2r+1)$-fold  sums of odd powers of the
first $n$ positive integers.

Applying Theorem \ref{sumrfacto} together with the following
relation
\begin{eqnarray}\label{eqn7}
(x+n)(x+n+k-1)_{2k-1}={1\over 2}(x+n+k)_{2k}+{1\over
2}(x+n+k-1)_{2k},
\end{eqnarray}
we arrive at the following formula for the $(2r)$-fold summation
 of even powers.

\begin{theorem}\label{thm3_2}
For $m\geq 1$, we have,
\begin{eqnarray*}
{\sum}^{2r}(x+n)^{2m}&=&\sum^m_{k=1}T(2m,2k)\Bigg\{\frac{(x+n+r)(x+n+k+2r-1)_{2k+2r-1}}{(2k+2r)_{2r}}\nonumber\\
&&-\sum^{2r}_{i=1}{n+2r-i-1\choose
2r-i}\frac{(2x+i)(x+k+i-1)_{2k+i-1}}{2(2k+i)_i}\Bigg\}.
\end{eqnarray*}
\end{theorem}

Setting $x=0$ in the above formula, the second summation vanishes.
One sees that the $(2r)$-fold summation becomes a polynomial in
$n(n+2r)$ because $(n+r)(n+k+2r-1)_{2k+2r-1}$ can be rewritten as
\[ \prod^{k+r}_{i=1}(n+k+2r-i)(n-k+i)\\[6pt]
=\prod^{k+r}_{i=1}[n(n+2r)+(k+2r-i)(i-k)].
\]

\section{$r$-Fold Alternating Sums of Powers}

In this section, we investigate the $r$-fold alternating sums of
powers.  Following the notation of Faulhaber, we define
\begin{equation}
{\sum}^{r}(-1)^n(n+x)^{m}:=\sum_{1\leq i_1 \leq \cdots \leq i_r \leq
n}(-1)^{i_1}(i_1+x)^{m}.
\end{equation}
We will  show that the $2r$-fold alternating sum of even powers
$\sum^{2r}(-1)^nn^{2m}$ is a polynomial in $n(n+r)$. For other cases
concerning the parities of $r, m$, we will outline the results
without proofs.

Define
\begin{equation}\label{er}
E^{(r)}_{m}(x_1,\ldots, x_{r}):=\sum_{i_1+\cdots +i_r=m}{m\choose
i_1,\dots, i_r }E_{i_1}(x_1)\cdots E_{i_r}(x_r).
\end{equation}

The following lemma holds. The proof is based on induction and is
omitted.

\begin{lemma}\label{lem-r-fold}
Let $r$, $m$ be positive integers. Then
\begin{eqnarray}
&&{\sum}^{r}(-1)^n(n+x)^{m}\nonumber\\
&&\quad=(-1)^n2^{-r}E^{(r)}_{m}\Big(1/2,\ldots,1/2,(x+n+r/2+1/2)\Big)\nonumber\\
&&\qquad +\sum^{r}_{k=1}{n+r-k-1\choose
r-k}2^{-k}E^{(k)}_{m}\Big(1/2,\ldots, 1/2,x+(k+1)/2\Big).\nonumber
\end{eqnarray}
\end{lemma}

We now give   recursive formulas for $E^{(k)}_{2m}$ in order to
compute the multiple sums in the above lemma.

\begin{lemma}\label{e-k-m}Let $k$, $m$ be positive integers.
Then $E^{(k)}_{2m}\left[1/2,\ldots,1/2,{(k+1)/ 2}\right]$ equals
\begin{equation}\label{t-2m}
\sum^{{k/2}}_{i=0}{k\choose 2i}\sum^i_{j=0}{i\choose
j}(-1)^jE^{(2j)}_{2m}\left(1/2,\ldots,1/2\right),
\end{equation}
and $E^{(k)}_{2m+1}\left(1/2,\ldots,1/2,{(k+1)/ 2}\right)$ equals
\begin{equation}\label{t-2m+1}
\sum^{{k/2}}_{i=0}{k\choose 2i+1}\sum^i_{j=0}{i\choose
j}(-1)^{j+1}E^{(2j+1)}_{2m+1}\left(1/2,\ldots,1/2,1\right).
\end{equation}
\end{lemma}

\pf From the generating function of $E_{m}(x)$, one sees that the
exponential generating function of $E^{(k)}_{m}\Big(1/2,\ldots,
1/2,(k+1)/2\Big)$ equals
\[
\sum^\infty_{m=0}E^{(k)}_{m}\Big(1/2,\ldots,
1/2,(k+1)/2\Big)\frac{t^m}{m!}=\left(\frac{2e^t}{1+e^t}\right)^k.
\]
Observe that
\[
\left(\frac{2e^t}{1+e^t}\right)^2=2\cdot\frac{2e^t}{1+e^t}-\left(\frac{2e^{{t\over
2}}}{1+e^t}\right)^2.
\]
Denote the exponential generating functions of $E_n(1)$ and
$E_n(1/2)$ by  $A(t)$ and $B(t)$, respectively. Then the above
relation implies that $A(t)^2=2A(t)-B(t)^2$, which yields
\[ A(t)=1-\sqrt{1-B(t)^2}.\]
Taking the $k$-th power,   we obtain
\begin{eqnarray*}
A(t)^k & = & \sum^{{k/2}}_{i=0}{k\choose 2i}\sum^i_{j=0}{i\choose
j}(-1)^jB(t)^{2j} \\[6pt]
& & \quad +\, \sum^{{(k-1)/2}}_{i=0}{k\choose
2i+1}\sum^i_{j=0}{i\choose j}(-1)^jB(t)^{2j}(1-A(t)).
\end{eqnarray*}

Since $E_{2m+1}(1/2)=0$ and $E_{2m}=0$ for $m\geq 1$, equating the
coefficients of both sides of above identity yields the desired
recurrence relations.
 \qed

The following theorem is concerned with  $2r$-fold alternating sums
of even powers.

\begin{theorem}
Let $r$, $m$ be positive integers.  Then the $2r$-fold alternating
sum $ \sum^{2r}(-1)^nn^{2m}$ is of the form
$(-1)^nF(n(n+2r))+G(n(n+2r))$, where $F$ and $G$ are polynomials of
degree $m$ and $r-1$ respectively.
\end{theorem}

\pf Applying Lemma  \ref{lem-r-fold} with $2r$ replaced by $2$, $m$
replaced by $2m$, and setting $x=0$,  we get
\begin{eqnarray}
&&{\sum}^{2r}(-1)^nn^{2m}\nonumber\\
&&\quad =(-1)^n2^{-2r}E^{(2r)}_{2m}
\Big(1/2,\ldots,1/2,(n+r+1/2)\Big)\label{2r-alt1}\\
&&\qquad \quad +\sum^{2r}_{k=1}{n+2r-k-1\choose
2r-k}2^{-k}E^{(k)}_{2m}\Big(1/2,\ldots,
1/2,(k+1)/2\Big).\label{2r-alt}
\end{eqnarray}
In the expansion of (\ref{2r-alt1}) indexed by
$2i_1+2i_2+\cdots+2i_{2r}=2m$, each term  contains a factor of the
form $E_{2i_r}(n+r+1/2)$. According to (\ref{e-bi}), we find that
\[
E_{2i_{2r}}\left(n+r+{1\over 2}\right)=
\sum^{i_{2r}}_{k=0}{2i_{2r}\choose 2k}E_{2i_{2r}-2k}\left({1\over
2}\right)(n+r)^{2k}.
\]
Since $(n+r)^{2k}=(n(n+2r)+r^2)^k$, it is easily seen that
 (\ref{2r-alt1}) is  a polynomial in $n(n+2r)$ of
degree $m$.
 We need to show that (\ref{2r-alt})
is also a polynomial in $n(n+2r)$. Applying (\ref{t-2m}), we find
that (\ref{2r-alt}) equals
\begin{eqnarray}
&&\sum^{2r}_{k=1}{n+2r-k-1\choose
2r-k}2^{-k}\sum^{{k/2}}_{i=0}{k\choose 2i}\sum^i_{j=0}{i\choose
j}(-1)^jE^{(2j)}_{2m}\left(1/2,\ldots,1/2\right)\nonumber\\
&&
\quad=\sum^{r}_{j=0}(-1)^jE^{(2j)}_{2m}\left(1/2,\ldots,1/2\right)\nonumber\\
&& \qquad \times \sum^r_{k=j}{n+2r-k-1\choose
2r-k}2^{-k}\sum^{k/2}_{i=j}{k\choose 2i}{i\choose j} .\label{kj}
\end{eqnarray}
Using the identity
\begin{equation}\label{bi-2i}
\sum^{n/2}_{i=j}{n\choose 2i}{i\choose j}=2^{n-2j-1}{n-j\choose
j}{n\over n-j},
\end{equation}
we deduce that the second sum in (\ref{kj}) equals
\begin{equation}\label{rkj}
\sum^r_{k=j}{n+2r-k-1\choose 2r-k}{k-j-1\choose j-1}{k\over
j}2^{-1-2j}.
\end{equation}
 In view of the following relation,
 \[
\sum_k{n-k\choose i}{m+k\choose j}={m+n+1\choose i+j+1}.
\]
the above sum (\ref{rkj})
 simplifies to
\begin{equation}\label{n2r}
 {n+2r-j-1\choose
2r-2j-1}\frac{n+r}{r-j}2^{-1-2j}.
\end{equation}
Now,  $(n+r)(n+2r-j-1)_{2r-2j-1}$ can be rewritten as
\[
\prod_{i=1}^{r-j}
(n+i+j)(n+2r-i-j)=\prod^{r-j}_{i=1}[n(n+2r)+(2r-j-i)(i+j)].
\]
which is a polynomial in $n(n+2r)$ of degree $r-j$. Since
$E^{(0)}_{2m}(\cdot)=0$, (\ref{2r-alt}) is a polynomial in $n(n+2r)$
of degree $r-1$. This completes the proof. \qed

For the remaining three cases with respect to the parities of $r$
and $m$, we have the following theorem. The proof is omitted.

\begin{theorem} For $r, m\geq 0$, we have
\begin{eqnarray}
&&{\sum}^{2r+1} (-1)^n
n^{2m}\label{2r+1-2m}\\
&&\quad=(-1)^nF^{(1)}_m
(n(n+2r+1))+(2n+2r+1)G^{(1)}_{r-1}(n(n+2r+1)),\nonumber\\[6pt]
&&{\sum}^{2r} (-1)^n n^{2m+1}\label{2r-2m+1}\\
&&\quad=(-1)^n(n+r)F^{(2)}_m(n(n+2r))+(n+r)G^{(2)}_{r-1}(n(n+2r)),\nonumber\\[6pt]
&&{\sum}^{2r+1} (-1)^n
n^{2m+1}\label{2r+1-2m+1}\\
&&\quad=(-1)^n(n+r)F^{(3)}_m(n(n+2r+1))+G^{(3)}_{r}(n(n+2r+1)),\nonumber
\end{eqnarray}
where $F^{(i)}_m(x)$ and $G^{(i)}_r(x)$ $(i=1,2,3)$ stand for
polynomials of degrees $m$ and $r$ respectively, and
$G^{(i)}_{-1}=0$ for $i=1,2$.
\end{theorem}

\vskip 3mm

\noindent \textbf{Acknowledgments.} We are grateful to the referees
for helpful comments leading to  improvements of an earlier version.
This work was supported by the 973 Project, the PCSIRT Project of
the Ministry of Education, the National Science Foundation,  and the
Ministry of Science and Technology of China.


\end{document}